\documentclass[11pt]{amsart}
\usepackage{amsmath, amscd, amssymb, latexsym}

\theoremstyle{plain}
\newtheorem{thm}{Theorem}[section]

\newtheorem{cor}[thm]{Corollary}
\newtheorem{pro}[thm]{Proposition}

\theoremstyle{definition}
\newtheorem{df}[thm]{Definition}

\newtheorem{rem}[thm]{Remark}

\def\om{\omega}
\def\Om{\Omega}

\def\ov{\overline}

\def\al{\alpha}

\def\si{\sigma}
\def\Si{\Sigma}
\def\na{\nabla}

\def\wt{\widetilde}

\def\ni{\noindent}
\def\nt{\noindent}

\def\un{\underline}

\def\pa{\partial}

\def\ti{\tilde}
\def\mc{\mathcal}

\def\mr{\mathrm}
\def\R{{\bf R}}
\def\C{{\bf C}}

\begin{document}

\vspace{-20mm}

\noindent{\em Comm. Contemp. Math.} {\bf 11} (2009) 109-130

\vspace{15mm}

\title{GERBES, HOLONOMY FORMS AND REAL STRUCTURES}

\author{SHUGUANG WANG}
\address{Department of Mathematics
\\University of Missouri\\Columbia,
MO 65211,USA}

\email{wangs@missouri.edu}

 \thanks{Work was partially supported by University of Missouri
Research Board Grant.}

\begin{abstract}
We study geometry on real gerbes in the spirit of
Cheeger-Simons theory. The concepts of  adaptations and holonomy
forms are introduced for flat connections on real gerbes. Their
relations to complex gerbes with connections are presented, as
well as results in loop and map spaces.
\end{abstract}

\maketitle

\nt{\em Keywords}: Real gerbes; adaptations; 
holonomy form; Cheeger-Simons; map space. 


\section{Introduction}\label{intro}

The concept of (complex) gerbes, especially non-abelian ones, was
first introduced by J. Giraud in 1971. Then the Chern-Weil theory
of gerbes was carried out in the book of Brylinski \cite{b}.
Abelian gerbes were further refined and clarified by Murray
\cite{m} and Chatterjee \cite{c} from two different view-points.
The two approaches  have their own advantages, and to some degree,
it is analogous to the bundle theory where one could opt for either
the principal or vector bundle approaches. Gerbes have found
applications in Physics for  the so-called $B$-fields and in
generalized index theory for twisted vector bundles,
see for example \cite{a, bm, mms} among many others.

In this paper, we will be undertaking the real version of
gerbes by adopting an approach more in line with Chatterjee
and Hitchin \cite{h}.
Just like real vector bundles, one should substitute the
Cheeger-Simons theory for the Chern-Weil theory. Thus we focus on
flat connections on real gerbes and introduce the concept of
adaptations in order to characterize the holonomy form of flat
connections. It turns out that one can relate real gerbes with
complex gerbes, provided the latter admit real structures in a
certain sense. This will be the second topic in the paper. Our
third topic is to extend the main results in Brylinski \cite{b}
and Chatterjee \cite{c} to the real gerbe case. These involve the
loop space of a manifold and the map space between a surface and a
manifold.

Here is an outline of the paper: in Section \ref{line}, we review
some basic facts on real line bundles with flat connections.  We
show how to get a flat connection from a holomorphic structure on
a complex line bundle (Theorem \ref{fla}). This will be used in
later sections. In Section \ref{ger}, we present the basic theory
on real gerbes and flat connections. The concepts of
adaptations and holonomy forms are introduced in Definition
\ref{gc}, as well as their important properties (Theorems
\ref{cur}, \ref{obj}). Similar to the case of line bundles, we
relate real and complex gerbes. The last section deals with the
loop and map spaces. By using the holonomy of a flat gerbe
connection, one constructs a line bundle on the loop space
(Theorem \ref{geo}). Through a certain boundary map, one pulls
back this bundle to the map space and the interesting result here
is that the pull-back line bundle is canonically trivial (Theorem
\ref{hl}). In the complex gerbe case of \cite{b}, such a result
has the origin in Topological Quantum Field Theory.

It is noted that throughout the paper, Riemannian metrics on
a real line bundle $l$ occupy a significant role, one reason
being that they lead to  canonical trivializations on
$l\otimes l$.

\section{Holonomy forms, real and complex line bundles}\label{line}

In this section, we start with a few basic but perhaps not exactly
familiar facts on real line bundles. We then relate real to
complex line bundles, which serves as a prototype for gerbes. Some of
the results are new and have not appeared elsewhere to our knowledge.

Let $l\to X$ be a real line bundle with structure group $\R^*$
over a smooth manifold. We emphasize that throughout the paper,
$X$ is allowed to be non-orientable or orientable but not
oriented. The space of connections on $l$ is affine modeled on
$\Om^1(X)$, since $\mr{End}l=l\otimes l^*$ is canonically trivial,
and acted by the gauge group $\mr{Map}(X,\R^*)$. Whereas the space
of flat connections on $l$ is affine modeled on the space $Z^1(X)$
of closed 1-forms and is preserved by the gauge group action.

\begin{pro}\label{cp}
For any fiber metric $h=(\;\;)$ on $l$, there is a unique
compatible connection $\nabla_h$, namely such that
$$d(s,t)=(\nabla_h s, t)+(s,\nabla_h t)$$
where $s,t$ are sections. Moreover $\nabla_h$ is flat.
\end{pro}

\ni{\it Proof}. Locally under a frame $s_i$ of $l$ with
$h_i=(s_i,s_i)>0$, the connection is given as
$\nabla_h=d+\frac{1}{2}d\ln h_i$, which is flat since
$\frac{1}{2}d\ln h_i$ is closed. \qed

In the case of a Hermitian complex line bundle, one usually 
concentrates
on unitary connections. Proposition \ref{cp} indicates that it
would be too restrictive to consider only Riemannian connections
in the real case. Instead one should expand to consider all flat
connections. Then one is led to measure how far a flat connection
is from being a metric one.

\begin{df}\label{cfm}
Given a flat connection $\nabla$ on $l$, the
{\em holonomy form} (or
{\em compatibility form}), with respect to a metric $h$,
 is defined to be the closed 1-form
$\theta(\nabla, h)=\nabla-\nabla_h\in Z^1(X)$. Locally under a
frame $s_i$,
\begin{equation}\label{lo}\theta(\nabla, h)=a_i-\frac{1}{2}d\ln h_i
\end{equation}
where one writes $\na=d+a_i$ under the same frame.
\end{df}

Note that even if $\na$ is not flat, $\theta(\nabla, h)$ is still
a well-defined 1-form, although it is not closed anymore.

More generally for a higher rank Riemannian vector bundle with a
flat connection, by taking determinant, one can define a similar
closed 1-form, which was  used in a fundamental way by
Bismut-Zhang \cite{bz} in connection with analytic torsions.

Of course $\theta(\nabla, h)$ will be trivial if $\nabla$ is 
compatible with $h$. In general one sees readily that the 
class $[\theta(\nabla, h)]\in H^1(X,\R)$
is independent of $h$, so it makes sense to denote it by
$\ti{c}(\na)$. One can think of $\ti{c}(\na)$ as an obstruction to
the existence of a metric compatible with $\na$. Note that
$\ti{c}(\na)$ is invariant under the gauge transformation hence it
descends to a map on the quotient of flat connections
under the gauge action,
$$\ti{c}:\mc{F}/\mc{G}\to H^1(X,\R).$$
This map is in fact bijective. Moreover the image of $\mc{H}\to
\mc{F}/\mc{G}, h\mapsto[\na_h]$ is precisely $\ti{c}^{-1}(0)$,
where $\mc{H}$ is the space of metrics on $l$. In particular there
exists a flat connection not compatible with any metric, provided
$H^1(X,\R)\not=\{0\}$.

Recall that to each flat connection $\na$, one can associate  its
holonomy class  $c(\na)\in H^1(X,\R^*)$ by using parallel
transport. The class $c(\na)$ relates to $\theta(\na,h)$ as
follows, explaining the terminology used in Definition \ref{cfm}.

\begin{pro}\label{rel}
{\em (a)} The homomorphism $\ln|\bullet|: \R^*\to \R$ induces one
from $H^1(X,\R^*)$ to $H^1(X,\R)$, under which the holonomy class
$ c (\na)$ is mapped to $\ti{c}(\na)$.

{\em (b)} The sign map $\R^*\to {\bf Z}_2=\{\pm1\}$ gives rise to
the homomorphism $H^1(X,\R^*)\to H^1(X,{\bf Z}_2)$, under which
$c(\na)$ is mapped to $w_1(l)$ hence independent of the flat
connection $\na$.

{\em (c)} Conversely the natural homomorphisms ${\bf
Z}_2\hookrightarrow \R^*, \mr{exp}:\R\to \R^*$ lead to the
decomposition
$$H^1(X,\R^*)=H^1(X,{\bf Z}_2)\times H^1(X,\R)$$
as Abelian groups, under which $ c (\na)=(w_1({l}),\ti{c}(\na))$.
\end{pro}

\nt{\em Proof}. All three parts amount essentially to the
following basic fact: along a loop $\rho: S^1\to X$, the holonomy
of $\na$ (not necessarily flat) is given by
\begin{equation}\label{rh}
\mr{hol}_\rho(\na)=\varepsilon\cdot\mr{exp}\int_{S^1}
\rho^*\theta(\na,h)
\end{equation}
where $\varepsilon=w_1({l})[\rho]$, $[\rho]\in H_1(X)$ and
$w_1(l)$ is viewed as a map $H_1(X)\to\{\pm 1\}$. Compare with the
complex line bundle case in \cite{mp}. Indeed since $\rho^*\na$ is
flat, there is a cover $\{I_i, 1\leq i\leq n\}$ of $S^1$ by open
intervals such that on each $I_i$, the pull-back bundle $\rho^* l$
is trivialized by a flat frame $s_i$ of $\rho^*\na$. (Actually one
can choose $n=2$.) Then under $s_i$, the connection 1-form of
$\rho^*\na$ is $a_i=0$ and
$$\rho^*\theta(\na, h)=\theta(\rho^*\na, \rho^*h)=
-\frac{1}{2}d\ln h_i$$
where we set $h_i=\rho^*h(s_i,s_i)$ which is subject to
$h_{i+1}/h_{i}=g^2_{i, i+1}$ on $I_i\cap I_{i+1}$, with $g_{i,
i+1}$ denoting the transition function of $\rho^* l$. Pick any
point $t_i\in I_i\cap I_{i+1}$ for $i=1,\cdots, n-1$ and $t_n\in
I_n\cap I_1$. For convenience, set $t_0=t_n, h_{n+1}=h_1$. One
calculates the right side of (\ref{rh}) as
\[\begin{array}{ll}
\hspace{-10mm}\varepsilon\cdot\mr{exp}\int_{S^1}
\rho^*\theta(\na,h)&= \varepsilon\cdot\mr{exp}
\sum^{n}_{i=1}\int^{t_i}_{t_{i-1}}-\frac{1}{2}d\ln h_i\\
&=\varepsilon\cdot\mr{exp}\sum_i\frac{1}{2}
\ln [h_i(t_{i-1})/h_i(t_{i})]\\
&=\varepsilon\cdot \prod_i\sqrt{h_i(t_{i-1})/h_i(t_{i})}\\
&=\varepsilon\cdot \prod^n_{i=1}\sqrt{g^2_{i, i+1}(t_i)
 h_i(t_{i-1})/h_{i+1}(t_{i})}\\
&=\varepsilon\cdot\prod^n_{i=1}[\varepsilon_i
g_{i, i+1}(t_i)]\cdot\sqrt{h_1(t_{0})/h_{n+1}(t_{n})}\\
&=\varepsilon\cdot\prod_{i}\varepsilon_i\cdot
\prod_i g_{i, i+1}(t_i)\\
&=\prod_i g_{i, i+1}(t_i),
\end{array}\]
where $\varepsilon_i=\mr{sgn}g_{i, i+1}(t_i)$. The last term is
the holonomy of $\na$ along $\rho$ by definition, thus verifying
equation (\ref{rh}). The signs cancel out because
$\varepsilon=w_1({l})[\rho]$ is the total sign change
$\prod^n_{i=1}[\mr{sgn}g_{i, i+1}(t_i)]$ along $\rho$.

Of course $c(\na)[\rho]=\mr{hol}_\rho(\na)$; the rest of the proof
is now pretty clear.\qed

Part (c) says that $w_1({l}), \ti{c}(\na)$ are complementary each
other in the sense that they capture respectively the topological
and geometrical aspects of the bundle  $l$. The picture here
contrasts with the Chern-Weil theory, where the real Chern classes
alone capture both the topological and the geometrical aspects of
a complex bundle. Of course the picture here should be viewed
as an instance of the Cheeger-Simons theory \cite{cs}.

It is also possible to describe the above results in terms of Cech
cohomology. As usual let $\un{\R}, \un{\R}^*$ denote the sheaves
of $\R^*$ or $\R$-valued functions on $X$, while $\R, \R^*$ denote
the constant sheaves. On any open cover of $X$, the transition
functions of $l$ is a co-closed 1-Cech cycle in $\un{\R}^*$ hence
determines a class $[l]\in\check{H}^1(X,\un{\R}^*)$. The short
exact sequence ${\bf Z}_2\to\un{\R}^*\to\un{\R}$ yields a long
exact sequence of cohomology groups, which in turn proves that
$\check{H}^1(X,\un{\R}^*)$ is  naturally isomorphic to $H^1(X,{\bf
Z}_2)$, because the fineness of $\un{\R}$ implies
$\check{H}^k(X,\un{\R})$ vanishes for $k>0$. Then under this
isomorphism, $[l]$ is mapped to $w_1(l)$. Likewise, by using local
$\na$-flat trivializations of $l$, the transitions  are local
constants and form a co-closed 1-Cech cycle in $\R^*$,
representing a class  $\check{c}(\na)\in\check{H}^1(X,\R^*)$.
Under the natural isomorphism
$\check{H}^1(X,\R^*)\to{H}^1(X,\R^*)$, $\check{c}(\na)$ becomes
the holonomy class ${c}(\na)$.

\vspace{3mm}

We now compare this real picture with the complex and holomorphic
pictures. First let $M$ be any smooth manifold with an involution
$\si$ and $L\to M$ a complex line bundle. Since $\si$ is to be
considered as a real structure, we assume $L$ admits a fiberwise
{\em conjugate} linear lifting $\tau$ of $\si$. We adopt the usual
convention to call $(L, \tau)$ a real complex line bundle or
simply a Real line bundle. It will turn out to be useful also to
take the conjugate linear lifting $\si^*: T^*M\otimes\C\to
T^*M\otimes\C$ on the $\C$-factor, hence a conjugate linear
extension $\si^*:\Om^p_\C(M)\to \Om^p_\C(M)$ on complex valued
forms. Then we  call  a $\si^*$-invariant complex form on $M$ a
{\em Real form} (as oppose to a real form in $\Omega^p(X)$). In
the same spirit, a $\tau$-invariant section of $L$ is refereed as
a {\em Real section} and in general, a Real object is simply a
complex object invariant under some real structure.

Set $X=M_\R:=\mr{Fix}\si$, which we assume to be a smooth
manifold. Set also the real line bundle $l=L_\R:=\mr{Fix}\tau\to
X$. Given a Real connection $d_A$ on $L$, namely $d_{\ov{A}}=d_A$
where by definition $d_{\ov{A}}s=\tau d_A\si^{-1}s$ for a section
$s$, it obviously restricts to a connection $\na_A$ on $l$. The
curvature 2-form $F$ of $d_A$ satisfies $\ov{F}=F$, where
$\ov{F}=\si^*F$. Thus ${F}$ restricts to a real 2-form $F'$ on $X$
($F'$ of course  is the curvature  of $\na_A$). If $d_A$ is also
unitary then $F$ is purely imaginary,  hence $F'$ must be trivial
and $\na_A$ is flat. Alternatively if $d_A$ is compatible with a
Real Hermitian metric $H$ on $L$ (i.e. $H(\ov{u},\ov{v})=H(u,v)$),
then $\na_A$  is compatible with $h$ hence flat, where the
Riemannian metric $h$  on $l$ is the restriction of $H$.

Now move on to the holomorphic picture and assume $M$ is a complex
manifold such that $\si$ is anti-holomorphic. A holomorphic
structure $\al$ on $L$ is characterized by its Dolbeault operator
$\ov{\pa}_\al:\Om^0(L)\to\Om^{0,1}(L)$, which satisfies
$$\ov{\pa}_\al(fs)=\ov{\pa}f\cdot s+f\cdot\ov{\pa}_\al s$$
and $\ov{\pa}_\al\circ \ov{\pa}_\al=0.$ The set of such operators
is affine modeled on the space $Z^{0,1}(M)$ of $\ov{\pa}$-closed
$(0,1)$ forms. This is quite analogous to flat connections on $l$.
Given a Hermitian metric $H$ on $L$, one is perhaps also tempted
to consider a compatibility  form $\Theta(\ov{\pa}_\al,H)\in
Z^{0,1}(M)$. Borrowing from formula (\ref{lo}), one might try to
set locally
$$\Theta(\ov{\pa}_\al,H)=\al_s-\frac{1}{2}\ov{\pa}\ln H_s$$
 where $\ov{\pa}_\al=\ov{\pa}+\al_s$ and
$H_s=H(s,s)>0$ under a local frame $s$ of $L$. (In particular, if
$s$ is holomorphic, $\Theta(\ov{\pa}_\al,H)=
-\frac{1}{2}\ov{\pa}\ln H_s$. Compare with the $\al, H$-compatible
connection $d_{\al,H}$, which is given by the local 1-form $\pa\ln
H_s$.) Of course $\Theta(\ov{\pa}_\al,H)$ is not going to be
well-defined, as the above expression depends on the choice of
$s$.  Nonetheless if $s$ is $\tau$-invariant,
$\Theta(\ov{\pa}_\al,H)$ can be restricted to a well-defined real
$1$-form on $X$, using which we prove the following result which
is somewhat surprising at first glance.

\begin{thm}\label{fla}
Any  holomorphic structure $\al$ on $L$ determines a unique flat
connection $\na_\al$ on $l$. In particular
 when $L$ is the trivial bundle with the trivial
$\al$,  $\na_\al$ recovers the differential $d$  on $X$.
\end{thm}

\ni{\it Proof.} We can assume $\al$ is Real, namely $\ov{\pa}_\al$
commutes with $\tau$; otherwise we replace $\ov{\pa}_\al$ by its
average $\ov{\pa}_\al+\frac{1}{2}(\tau^*(\ov{\pa}_\al)-
\ov{\pa}_\al)$. Take any Real Hermitian metric $H$ on $L$. Let $h$
be its restriction to $l$ and $\na_h$ the compatible flat
connection. The idea of proof is to first show that under a local
$\tau$-invariant frame $s$ of $L$ near $X$, the $(0,1)$ form
$$\Theta_s(\ov{\pa}_\al,H):=\al_s-\frac{1}{2}\ov{\pa}\ln H_s$$
restricts to a well-defined real closed 1-form
$\theta(\ov{\pa}_\al,H)$ on $X$, namely independent of $s$. Then
show that the flat connection defined as
\begin{equation}\label{tt}
\na_\al:=\na_h+{2}\theta(\ov{\pa}_\al,H)
\end{equation}
is independent of $H$.

Since $\al, H$ are both Real, $\al_s$ and $\ov{\pa}\ln H_s$ are
Real;hence $\Theta_s(\ov{\pa}_\al,H)$ restricts to a real local
1-form $\theta_s(\ov{\pa}_\al,H)$  on $X$. Let $s'=fs$ be another
$\tau$-invariant frame, so $H_{s'}=|f|^2 H$. Then
$$\begin{array}{ll}
\Theta_{s'}(\ov{\pa}_\al,H)&=\al_{s'}-\frac{1}{2}\ov{\pa}\ln H_{s'}\\
&=\al_s+f^{-1}\ov{\pa}f-\frac{1}{2}\ov{\pa}\ln H_{s}-\ov{\pa}\ln|f|\\
&=\Theta_{s}(\ov{\pa}_\al,H)+f^{-1}\ov{\pa}f-|f|^{-1}\ov{\pa}|f|.
\end{array}$$
When restricted to $X$, $f$ is real and positive.
Hence$f^{-1}\ov{\pa}f-|f|^{-1}\ov{\pa}|f|=0$, from which it
follows that
$\theta_{s'}(\ov{\pa}_\al,H)=\theta_s(\ov{\pa}_\al,H)$, giving a
well-defined form $\theta(\ov{\pa}_\al,H)$ on $X$. The form
$\theta(\ov{\pa}_\al,H)$ is closed since $\al_s $ and $\ov{\pa}\ln
H_s$ are $\ov{\pa}$-closed as well as Real.

Next we show that $\na_\al$ as defined in (\ref{tt}) is
independent of $H$. Let $H'=e^{f}H$ be a second Real Hermitian
metric where $f$ is some real valued $\si$-invariant function on
$M$. For their restrictions to $l$, $h'=e^{f}h$ holds as well. We
need to show
$$\na_h+{2}\theta(\ov{\pa}_\al,H)=\na_{h'}+
{2}\theta(\ov{\pa}_\al,H').$$ As before, $H'_s, h_s, h'_s$ denote
the norm squares of a local frame $s$ under the various metrics so
$H'_s=e^{f}H_s$ and $h'_s=e^{f}h_s$ are still valid. Then locally
we need $$\frac{1}{2}\ln h_s+{2}[\al_s- \frac{1}{2}\ov{\pa}\ln
H_s]|_X= \frac{1}{2}\ln h'_s+{2}[\al_s-\frac{1}{2}\ov{\pa}\ln
H'_s]|_X$$ which simplifies to $2\ov{\pa}f=df$ on $X$. The last
holds true since $\ov{\pa} f$ restricts to a real 1-form:
$\ov{\pa} f=\ov{\ov{\pa}{f}}={\pa}\ov{f}={\pa}{f}$ on $X$.

When $\al$ is trivial, we set $H_s=1$ under the global trivial
holomorphic frame so that $h_s=1, \theta(\ov{\pa}_\al,H)=0$ and
$\na_\al=d$. \qed

\begin{cor}
Let $\mc{F}^c$ denote the space of all holomorphic structures on
$L$ and $\mc{F}$ be the space of all flat connections on $\ell$.
The following diagram commutes:
$$\begin{array}{ccc}
\mc{F}^c\times Z^{0,1}(M)&\longrightarrow& \mc{F}^c\\
\downarrow&&\downarrow\\
\mc{F}\times Z^{1}(X)&\longrightarrow& \mc{F},
\end{array}$$
where the map $Z^{0,1}(M)\to Z^{1}(X)$ sends $\om$ to the
restriction of $(\om+\ov{\om})/2$ on $X$.
\end{cor}

It is possible to construct $\na_\al$ in Theorem \ref{fla}
directly without using a Hermitian metric. Let $s$ be
$\tau$-invariant local holomorphic frame  of $\al$. Then $s$
restricts to a local frame of $l$. One simply declares $t$ to be
the flat frame of a connection $\na_\al$ on $l$; one can check
that $\na_\al$ is well-defined, i.e. independent of the choice of
$s$. The flatness of $\na_\al$ follows from $\ov{\pa}_\al\circ
\ov{\pa}_\al=0$. Since $\al$ is Real, $\ov{\pa}_\al$ induces a map
between the spaces of Real forms: $\Om^0(L)_\R\to
\Om^{0,1}(L)_\R$. It is possible to show that the restriction maps
$\Om^0(L)_\R\to \Om^0(l), \Om^{0,1}(L)_\R\to \Om^1(l)$ are both
surjective. Then $\na_\al$ is the unique map making the following
diagram commute:
$$\begin{array}{ccc}&\ov{\pa}_\al&\\
\Om^0(L)_\R&\longrightarrow& \Om^{0,1}(L)_\R\\
\downarrow&&\downarrow\\
&\na_\al&\\
\Om^0(l)&\longrightarrow& \Om^{1}(l).
\end{array}$$

Clearly the theorem does not hold for a higher rank holomorphic
bundle $E$ with a real lifting, since the real bundle $E_\R$ may
not be flat in the first place.

The theorem suggests that the flat connection $\na_\al$  can be
considered appropriately as the counter-part of a holomorphic
structure $\al$. In the holomorphic set up, $\al$ and $ H$
determine a unique connection $d_{\al,H}$, while in the real set
up, $\na_\al$ and $h$ determine a unique 1-form $\theta(\na_\al,
h)$. (By equation (\ref{tt}), $\theta(\na_\al,
h)=2\theta(\ov{\pa}_\al,H)$. The flat connection
$\na_h+\theta(\ov{\pa}_\al,H)$ will depend on $H$. Note that
$d_{\al,H}$ restricts to the metric connection $\na_h$, hence
dependent on $H$ but having nothing to do with $\al$.)

Thus a holomorphic structure $\al$ on $L$ determines the
topological property of the line bundle $l$ via $\na_\al$. As an
application, one could study the link between the complex and real
analytic torsions developed by Bismut {\em et al} \cite{bgs,bl},
which are defined on holomorphic and real determinant line bundles
of elliptic operators. We hope to follow up this in a future work.

\section{Real gerbes, connections, and adaptations}\label{ger}

Real gerbes are next in the hierarchy after real line bundles. On
a smooth manifold $X$, a real gerbe $\mc{G}=\{U_i, l_{ij},
s_{ijk}\}$ consists of the following data:

\nt $\bullet$ $\{U_i\}$ is an open cover of $X$,

\nt $\bullet$ $l_{ij}$ is a real line bundle on $U_i\cap U_j$ with
a given isomorphism $l_{ij}\otimes l_{ji}=\wt{\R}$, 

\nt $\bullet$ $s_{ijk}$ is a trivialization of the bundle
$l_{ijk}=l_{ij}\otimes l_{jk}\otimes l_{ki}$ over $ U_i\cap
U_j\cap U_k$. Take a refinement cover if necessary so that each
$l_{ij}$ is a trivial bundle and $l_{ijk}$ has a second
trivialization. Then $s=\{s_{ijk}\}$ can be viewed as a Cech 2-cycle
of the sheaf $\un{\R}^*$. One imposes that $s$ be co-closed:
$\delta s=1$.

Thus the {\em gerbe class} $[\mc{G}]:=[s]\in\check{H}^2(X,
\un{\R}^*)$ is defined. We say that two gerbes are isomorphic if
they become the same on a common refinement. Then the isomorphic
class of $\mc{G}$ is determined by $[\mc{G}]$. Note that one can
also view $[\mc{G}]\in {H}^2(X, {\bf Z}_2)$ using the natural
isomorphism ${H}^2(X, {\bf Z}_2)\to\check{H}^2(X, \un{\R}^*)$,
which comes from the short exact sequence ${\bf
Z}_2\to\un{\R}^*\to \un{\R}$ as seen  before.\\

\nt{\em Examples}. (1) Any manifold $X$ is locally spin, hence
admits a collection of spinor bundles $S_i\to U_i$ which satisfy
$S_i=S_j\otimes l_{ij}$ for some line bundle $l_{ij}$ on $U_{ij}$.
In the language of \cite{ma} for example, $\{S_i\}$ form a twisted
vector bundle, the failure of which to be a global bundle is the
real gerbe $\{l_{ij}\}$ with the gerbe class $w_2(X)$.

(2) Take any codimension 2 submanifold $\Si\subset X$.
Trivializing the tubular neighborhood, one can construct easily a
real gerbe that has the gerbe class  PD$[\Si]\in H^2(X,{\bf
Z}_2)$.\\

So far this is in a complete analogy with a complex gerbe. On the
other hand, to do geometry we need a connection etc., which will
be different from the complex case. A gerbe connection
$\na=\{\na_{ij}\}$ on $\mc{G}$ consists of connections $\na_{ij}$
on $l_{ij}$ such that $s_{ijk}$ is a covariant constant
trivialization under the induced product connection
$\na_{ijk}=\na_{ij}\otimes\na_{jk}\otimes\na_{ki}$. We call $\na$
{\em flat} if all $\na_{ij}$ are flat. Likewise, a gerbe metric
$h=\{h_{ij}\}$ on $\mc{G}$ consists of a family of fiber metrics
$h_{ij}$ on $l_{ij}$ such that each trivialization $s_{ijk}$ has
norm 1 under the induced product metric $h_{ijk}$. The usual
argument of partitions of unity shows that connections and metrics
always exit on any gerbe $\mc{G}$. Furthermore, a metric $h$ leads
to a flat connection $\na_h=\{\na^h_{ij}\}$, where $\na^h_{ij}$ is
the flat connection compatible with $h_{ij}$. Hence we have
established the existence of flat gerbe connections on any gerbe
as well.

\begin{rem} In the complex case \cite{c,h}, a gerbe connection
$\na^c=\{\na^c_{ij},F_i\}$ consists of two parts: the part of
local connections $\na^c_{ij}$ on $U_i\cap U_j$ (called the
0-connection),  and the part of local 2-forms $F_i$ on $U_i$
(called the 1-connection) which are  subject to a certain
condition. Likewise it is possible to incorporate a collection of
local 2-forms into our real gerbe connection $\na=\{\na_{ij}\}$. 
However our main interest lies in flat gerbe connections, for 
which we can simply take trivial local 2-forms. Thus 1-connections 
do not play any significant role here.

\end{rem}

Given a flat gerbe connection $\na=\{\na_{ij}\}$,  $\na_{ij}$-flat
trivializations on $l_{ij}$  induce a $\na_{ijk}$-flat
trivialization on each $l_{ijk}$, which together with the original
$\na_{ijk}$-flat trivialization $s_{ijk}$ define a closed Cech
2-cycle of the constant sheaf $\R^*$. We call the resulted class
$c(\na)\in \check{H}^2(X, \R^*)=H^2(X,\R^*)$ the {\em holonomy} of
$\na$.  Given a smooth map defined on a closed oriented surface,
$f: \Si\to X$, the holonomy around $\Si$ is by definition the
value $c(\na)(f_*[\Si])\in \R^*$. As in the case of line bundles,
under the natural homomorphism $\check{H}^2(X,\R^*)\to
\check{H}^2(X,\un{\R}^*)$, the image of  $c(\na)$ recovers the
gerbe class $[\mc{G}]$. (Indeed both classes are represented by
the same 2-cocycle.)  Let $\ti{c}(\na)\in H^2(X,\R)$ be the image
class of $c(\na)$ under the other natural homomorphism
$H^2(X,\R^*)\to H^2(X,{\bf R})$. To do a Cheeger-Simons type
theory here means to find a differential form representing
$\ti{c}(\na)$. Such a form will be the gerbe version of the
holonomy form introduced in the previous section. The following is
a main definition of the paper.
\begin{df}\label{gc}
Consider a gerbe metric $h=\{h_{ij}\}$ and a flat connection
$\na=\{\na_{ij}\}$ on $\mc{G}$. An {\em adaptation} to $(\na, h)$
is a collection of 1-forms $\beta=\{\beta_i\} \mbox{ on }\{U_i\}$
such that
$$\beta_i-\beta_j=\theta(\na_{ij}, h_{ij})\;\;\mr{ on }\;\;
U_i\cap U_j,$$ where $\theta(\na_{ij}, h_{ij})$ is the holonomy
form of $\na_{ij}$ with respect to $h_{ij}$. Since each
$\theta(\na_{ij}, h_{ij})$ is closed, $\{d\beta_i\}$ fit together
to yield a  global 2-form $B$ on $X$. We call $B$ the {\em
holonomy form} of $(\na, h, \beta)$ (or of $\na$ with respect to
$h,\beta$).
\end{df}

Note that even if $\na$ is not flat, it still makes sense to
define adaptations. Of course the $2$-form $B$ is no longer
well-defined globally but can be viewed instead as a ``twisted''
2-form on $X$ over $\{d\theta(\na_{ij}, h_{ij})\}$.

The standard sheaf theory shows the existence of adaptations:
Since $\na, h$ are both compatible with the trivialization
$s_{ijk}$ locally, $\{\theta(\na_{ij}, h_{ij})\}$ is a closed Cech
1-cocycle in the sheaf $\mc{A}^1$ of 1-forms on $X$. The fineness
of $\mc{A}^1$ guarantees the existence of a 0-cycle
$\beta=\{\beta_i\}$ of 1-forms.

Adaptations are however not unique for a given metric and a flat
connection: any two  differ by a global 1-form on $X$. To some
degree, $\beta$ resembles a 1-connection in the complex gerbe
case, while $B$ is analogous to the curvature 3-form of the
1-connection. The following is parallel to Proposition \ref{rel}.

\begin{thm}\label{cur} Suppose $B$ is the holonomy form of
a flat gerbe connection $\na$ with respect to a gerbe metric $h$
and an adaptation $\beta$. Then the class $[B]\in H^2(X,\R)$
depends on $\na$ only. In fact $B$ represents the class
$\ti{c}(\na)$. Under the natural product $H^2(X,\R^*)=H^2(X,{\bf
Z}_2)\times H^2(X,\R)$, we  have $c(\na)=([\mc{G}],[B])$.
\end{thm}

\nt{\it Proof.} First we show that $[B]$ is independent of the
choice of adaptations: If $\beta'=\{\beta'_i\}$ is another
adaptation to $(\na, h)$, then there is a global 1-form $\al$ on
$X$ such that $\beta'_i=\beta_i+\al$. Then the corresponding
holonomy form $B'$ satisfies $[B']=[B+d\al]=[B]$.

To see $[B]$ is independent of $h$, let $h'=\{h'_{ij}\}$ be
another gerbe metric. Then there is a positive function $f_{ij}$
on each $U_i\cap U_j$ such that $h'_{ij}=f_{ij}h_{ij}$. Since the
metrics $h_{ijk}, h'_{ijk}$ both normalize the trivialization
$s_{ijk}$, $\{f_{ij}\}$ form a co-closed Cech 1-cycle of the sheaf
$\un{\R}^+$. But $\check{H}^1(X, \un{\R}^+)$ is trivial (due to
the fineness of $\un{\R}^+$), consequently, there is a 0-cycle
$\{f_i\}$ with co-boundary $\{f_{ij}\}$. Now one checks easily
that $\beta'_i=\beta_i-\frac{1}{2}d\ln f_i$ is an adaptation to
$(\na, h')$ and the associated holonomy form
$B'=d\beta'_i=d\beta_i=B$, from which $[B']=[B]$ for sure.

Let $\Si$ be a closed surface and  $f:\Si\to X$ a smooth map
representing a class $\al\in H_2(X, {\bf Z})$. The remaining
statement in the theorem is essentially due to that
 the holonomy of $\na$ around $\Si$ is
\begin{equation}\label{av}
c(\na)(\al)=\varepsilon\cdot \mr{exp}\int_\Si f^*B,
\end{equation}
where  $\varepsilon=[\mc{G}](\al)\in\{\pm1\}$. By functorality, it
is enough to show for the case that $X=\Si$ and $f=id$, where
(\ref{av}) becomes $c(\na)=\varepsilon\cdot \mr{exp}\int_\Si B$.
To compute the integral on the right side, the idea  is to
partition $\Si$ suitably and apply Stokes' Theorem repeatedly. For
example when $\Si$ is the 2-sphere $S^2$, view $S^2$ as a
rectangle $I^2$ with all four sides collapsed to the base point of
$S^2$. Take a small enough rectangular subdivision of $I^2$ so
that $\mc{G}$ is trivialized around all grid lines. On each
sub-rectangle, apply Stokes' Theorem twice to reduce the integral
first to the four sides and then to the four corners. The final
outcome is exactly the value $c(\na)\in
\check{H}^2(S^2,\R^*)=\R^*$. The complete details are left to the
interested reader (compare with the complex gerbe case in
\cite{mp}). \qed

\begin{cor}\label{hv}
Fix any metric $h$ on $\mc{G}$. A flat gerbe connection $\na$ on
$\mc{G}$ has  trivial holonomy   iff the gerbe class $[\mc{G}]$ is
trivial and there exists an adaptation  to $\na, h$ with the
holonomy form  vanishing identically.
\end{cor}

\nt{\em Proof}. In general, start with adaptation $\beta$ to
$(\na, h)$ with holonomy form $B$.  Then every form $B'$
representing the class $[B]$ can be realized by the holonomy form
of another adaptation $\beta'$ to $(\na, h)$. In fact,
$B'=B+d\gamma$ for some $1$-form $\gamma$ on $X$ and one simply
takes $\beta'=\{\beta'_i+\gamma\}$ where $\beta=\{\beta_i\}$ with
respect to some open cover $\{U_i\}$. In the situation of the
corollary, if $c(\na)=0$ then $[B]=0$ and one uses $B'=0$. The
rest is clear.\qed

The following sums up the main properties that will be quite
useful for Section \ref{hlo}.

\begin{thm}\label{obj}
Suppose the gerbe class $[\mc{G}]\in \check{H}^2(X,\un{\bf R}^*)$
is trivial.

{\em (a)} Then $\mc{G}$ admits a trivialization, namely a
collection of line bundles $\{l_i\}$ on $\{U_i\}$ together with
isomorphisms $l_i\otimes l^*_j=l_{ij}$ on $U_i\cap U_j$. Given a
second trivialization $\{l'_{i}\}$, a global line bundle $\xi$ on
$X$ is resulted by patching all $l_{i}\otimes(l'_{i})^*$ together.

{\em (b)} Given any gerbe connection $\{\na_{ij}\}$ on $\mc{G}$,
there is a collection of connections $\na_i$ on $l_i$ such that
$\na_i\otimes\na^*_j=\na_{ij}$ under the isomorphism $l_i\otimes
l^*_j=l_{ij}$. From a second trivialization $\{l'_i\}$ with
connections $\{\na'_i\}$, $\na_i\otimes(\na'_i)^*$ together form a
well-defined global connection $D$ on the bundle $\xi$.

{\em (c)} Suppose further that $\na=\{\na_{ij}\}$ is a flat gerbe
connection with trivial holonomy, $c(\na)=0$. Then one can choose
all connections $\{\na_i\}$ in part {(b)} to be flat. Given a
second trivialization with flat local connections, the induced
connection $D$ from the last part is flat as well.

{\em (d)} Let $h=\{h_{ij}\}$ be a gerbe metric on $\mc{G}$. Then
there is a collection of metrics $h_i$ on $l_i$ such
that$h_i/h_j=h_{ij}$ on $U_i \cap U_j$. For the flat connections
$\{\na_i\}$ constructed in part (c), the 1-form collection
$\beta^0=\{\theta(\na_i, h_i)\}$ is an adaptation to $\na$ and
$h$. Moreover the holonomy form $B^0$ of $(\na,h,\beta^0)$ is
trivial, which in particular proves Corollary \ref{hv} for a
second time.

{\em (e)} The second part of (d) has a partial converse in the
following sense. Fix a gerbe metric $h=\{h_{ij}\}$ and a
collection of local metrics $\{h_i\}$ as in (d). Suppose
$\na=\{\na_{ij}\}$ is a flat gerbe connection with trivial
holonomy and $\beta=\{\beta_i\}$ is any adaptation to $(\na, h)$.
If  holonomy form $B=$ of $(\na, h,\beta)$ is trivial, then there
exit a possibly different flat gerbe connection
$\na'=\{\na'_{ij}\}$ with trivial holonomy and a collection of
flat local connections $\{\na'_i\}$ subordinate to $\na'$ as in
(c) such that $\beta_i=\theta(\na'_i,h_i)$. In other words,
(d) and (e) imply essentially that for a given $\beta$, equations
$\beta_i=\theta(\na'_i,h_i)$ admit solutions for $\na'_i$ iff
$B=0$.
\end{thm}

\ni{\em Proof.} (a) Take a refinement cover of $\{U_i\}$ if
necessary, so that each $l_{ij}$ is trivialized on $U_{ij}$ and
$\mc{G}$ is represented by a 2-cocycle $s=\{s_{ijk}\}$. Since
$[\mc{G}]=[s]=0\in \check{H}^2(X,\un{\R}^*)$, there is a Cech
1-cocycle $f=\{f_{ij}\}\in \check{C}^1(X,\un{\R}^*)$ with
coboundary $\delta f=s$. Then  the trivial bundle $l_i\to U_i$ is
glued with the trivial bundle $l_j\otimes l_{ij}$ on $U_i\cap U_j$
via $f_{ij}$. By utilizing $f_{ij}$ as local transition functions,
this gives a desired trivialization on the original cover when
$l_{ij}$ is not necessarily trivial. For another trivialization
$\{l'_i\}$ of $\mc{G}$, we have $l_i\otimes l^*_j=l'_i\otimes
(l'_j)^*$ on $U_i\cap U_j$. Hence $l_i\otimes (l'_i)^*=l_j\otimes
(l'_j)^*$, namely the local bundles $\{l_i\otimes (l'_i)^*\}$ glue
together to form a global real line bundle $\xi$ on $X$.

(b) Still assume $\mc{G}=\{l_{ij}\}$ and $\{l_i\}$ are both
locally trivialized as in (a), so that $\na_{ij}=d+a_{ij}$ and
$\na_{i}=d+a_{i}$ for some 1-forms $a_{ij}, a_{i}$. To find the
required $\na_i$, one needs to have some $a_{i}$ such that
\begin{equation}\label{ac}
a_i-a_j=a_{ij}+f_{ij}^{-1}df_{ij},
\end{equation}
namely $T:=\{a_{ij}+f_{ij}^{-1}df_{ij}\}$ is a coboundary  Cech
cycle in the  sheaf $\mc{A}^1$ of 1-forms on $X$. This is the case
as $T$ is  closed  and the Cech cohomology
$\check{H}^1(X,\mc{A}^1)$ is trivial from the fineness of
$\mc{A}^1$. For a second trivialization $\{l'_i\}$ with local
connections $\{\na'_i\}$, clearly
$\na_i\otimes\na'_i=\na_j\otimes\na'_j$ under the same isomorphism
$l_i\otimes (l'_i)^*=l_j\otimes (l'_j)^*$, forming a connection
$D$ on $\xi$ by gluing.

(c) When $\na=\{\na_{ij}\}$ is flat, we choose flat
trivializations for each bundle $l_{ij}$ in parts (a), (b) above,
so that $a_{ij}=0$ and $\mc{G}$ is represented by the 2-cocycle
$s=\{s_{ijk}\}$ which now lives  in the constant subsheaf
$\R^*\subset \un{\R}^*$. Since the holonomy $c(\na)=[s]=0\in
\check{H}^2(X, \R^*)$, one can choose a local constant 1-cycle
$\{f_{ij}\}$ with coboundary equal to $s$. Thus in equation
(\ref{ac}) above, the right side is identically zero, which means
we can choose all $a_i=0$ to get the desired flat connections
$\{\na_i\}$. For a second trivialization of $\mc{G}$ together with
local flat connections, the induced connection $D$ is flat.

(d) As above, refine the open cover so that $\mc{G}$ is
trivialized locally. Then each $h_{ij}$ is a positive function on
$U_i\cap U_j$. Since $h_{ij}h_{jk}h_{ki}=1$ on $U_i\cap U_j\cap
U_k$, $\{h_{ij}\}$ form a closed 1-cocycle in the sheaf
$\un{\R}^+$ of positive functions. As $\check{H}^1(X, \un{\R}^+)$
vanishes by fineness of $\un{\R}^+$, there is a 0-cycle $\{h_i\}$
with coboundary $h_ih^{-1}_j=h_{ij}$.

For flat connections $\na_i$ constructed in (c), by definition
$\theta(\na_i, h_i)=a_i-\frac{1}{2}d\ln h_i$ and $\theta(\na_{ij},
h_{ij})=a_{ij}-\frac{1}{2}d\ln h_{ij}$. It is then easy to see
that $\{\theta(\na_i, h_i)\}$ is an adaptation to $\na$ and $ h$,
namely
$$\theta(\na_i, h_i)-\theta(\na_j, h_j)=\theta(\na_{ij},h_{ij}).$$
The holonomy form $B^0$ is trivial, since $\theta(\na_i, h_i)$ are
all closed.

(e) Start with some local flat connections $\{\na_i\}$ subordinate
to $\na$ as in (c). Set $\gamma_i=\beta_i-\theta(\na_i, h_i)$.
Then on $U_i\cap U_j$,
$$\begin{array}{ll}
\gamma_i-\gamma_j&=[\beta_i-\beta_j]-[\theta(\na_i, h_i)
-\theta(\na_j, h_j)]\\
&=\theta(\na_{ij}, h_{ij})-\theta(\na_{ij}, h_{ij})=0.
\end{array}$$
Hence we have a global 1-form $\gamma$ on $X$ and $d\gamma=B$.

If $B=0$, then on a possibly refined open cover,
$\beta_i-\theta(\na_i, h_i)=dg_i$ for some function $g_i$. Put
$g_{ij}=g_i-g_j$ on $U_i\cap U_j$, and define
$\na'_{ij}=\na_{ij}+dg_{ij}$. Since  $\{g_{ij}\}$ is a closed Cech
1-cycle in the sheaf $\un{\R}$, $\{\na'_{ij}\}$ form a flat gerbe
connection with trivial holonomy. Moreover the flat local
connections $\na'_i=\na_i+dg_i$  are subordinate to
$\{\na'_{ij}\}$ in the sense of part (c). One checks readily that
$\beta_i=\theta(\na'_i, h_i)$. \qed

In (c), without assuming $c(\na)=0$, there may not exist flat
local connections $\na_i$  satisfying (\ref{ac}), because the
sheaf $\mc{Z}^1$ of closed 1-forms on $X$ has the non-trivial
cohomology $\check{H}^1(X,\mc{Z}^1) $ in general. (In fact
$\check{H}^1(X,\mc{Z}^1)=H^1(X,\R)$.)

\begin{df}
Adapting the terms from Chatterjee \cite{c}, we call  $\{l_i\},
\{\na_i\}, \{h_i\}$ an {\em object bundle, object connection} and
{\em object metric}, which will all be referred to as objects
conveniently. They are respectively {\em subordinate to} the gerbe
$\mc{G}$, the connection $\na$ and the metric $h$. Furthermore,
two objects of $\mc{G}$ determine the {\em difference bundle}
$\xi$, two objects of $\na$ determine the {\em difference
connection} $D$, etc.
 \end{df}

From another point of view, as in \cite{ma,y} for example,
$\{l_i\}$ can also be called a twisted line bundle on $X$ over the
gerbe $\mc{G}$, while $\{\na_i\}$ a twisted connection  over the
gerbe connection $\na=\{\na_{ij}\}$.

In the final part of the section, we consider Real complex gerbes
and relate them to real ones. As in the previous section,
$\si:M\to M$ is a smooth involution with fixed point set $X$. Take
a complex gerbe $\mc{G}^c=\{L_{ij}, s^c_{ijk}, U^c_i; i,j, k\in
I\}$ on $M$. Assume the cover $\{U^c_i; i\in I\}$ is Real, namely
there is  an involution $I\to I, i\mapsto\ov{i}$, such that
$\si:U^c_i\to U^c_{\ov{i}}$ is a diffeomorphism for any $i\in I$,
and $i=\ov{i}$ whenever $X\cap U^c_{i}\not=\emptyset$. Then we say
$\mc{G}^c$ is Real if $\si: U^c_i\cap U^c_j\to U^c_{\ov{i}}\cap
U^c_{\ov{j}}$ has a Real lifting $\tau_{ij}: L_{ij}\to
L_{\ov{i}\:\ov{j}}$. Obviously the real part of $\mc{G}^c$ yields
a real gerbe $\mc{G}=\{l_{ij}, s_{ijk}, U_i; i,j,k\in I_\R\}$ on
$X$ where $I_\R=\{i\in I; X\cap U^c_{i}\not=\emptyset\}$.

A gerbe connection $\na^c=\{\na^c_{ij},F_i\}$ on $\mc{G}^c$ is
called Real if each $\na^c_{ij}$ is Real with respect to
$\tau_{ij}$ and $\ov{F_i}=F_i$ (i.e. $F_i$ is Real also), where $
\ov{F_i}={\si^*F_{\ov{i}}}$ and  $\si^*$ is conjugate linear on
complex forms as before. A Hermitian metric $H=\{H_{ij}\}$ on
$\mc{G}^c$ is called Real in the obvious sense.

According to Chatterjee \cite{c}, a gerbe connection
$\na^c=\{\na^c_{ij},F_i\}$ is compatible with a gerbe Hermitian
metric $H$ if each $\na^c_{ij}$ is compatible with $H_{ij}$ in the
usual sense and $F_i$ is purely imaginary, the latter being  so
required in view that the curvature of each $\na^c_{ij}$ is purely
imaginary. Clearly a Real Hermitian connection $\na^c$ restricts
to a flat connection $\na=\{\na_{ij}; i,j\in I_\R\}$ on $\mc{G}$.
Note that $F_i$ restricts trivially on $X$, suggesting once more
that there is no need to include any 2-forms as a part of any real
gerbe connection.

Suppose further that $M$ is a complex manifold. By definition a
holomorphic structure $\al=\{\al_{ij}\}$ on $\mc{G}^c$ consists of
a collection of holomorphic structures $\al_{ij}$ on $L_{ij}$ such
that $s^c_{ijk}$ is a holomorphic section. A connection $\na^c$ is
compatible with $\al$ if each $\na^c_{ij}$ is so with $\al_{ij}$
and $F_i$ has no $(0,2)$-component (see \cite{c}). Thus a
Hermitian holomorphic gerbe $\mc{G}^c$ admits a unique compatible
set $\{\na^c_{ij}\}$ of local connections (the 0-connection),
while the 1-connection $\{F_i\}$ is now a collection of local
imaginary $(1,1)$-forms which  however  are not unique.

Assume $\si: M\to M$ is anti-holomorphic and $\al$ is Real in the
sense that each $\al_{ij}$ is Real with respect to $\tau_{ij}$.
Then Theorem \ref{fla} says that any Real holomorphic structure
$\al$ on $\mc{G}^c$ determines a unique flat connection $\na^\al$
on $\mc{G}$. Note that $\al, H$  do not determine an adaptation
subordinated to $(\na^\al, h)$, nor does a 1-connection $\{F_i\}$.
However we can obtain a unique adaptation from objects:\\

\nt{\em Example}. Suppose the holomorphic gerbe class $[\al]=0\in
H^2(M,\mc{O}^*)$ so that $(\mc{G}^c, \al)$ admits a holomorphic
object $\{L_i, \al_i\}$. Endow $\{L_i\}$ with a Hermitian object
metric $H^{ob}=\{H_i\}$ of $H=\{H_{ij}\}$. If the objects are both
Real, then one  has a well-defined adaptation to $(\na^\al,h)$,
given by $\{\theta(\na^{ob}_i,h_i)\}$, where $\na^{ob}_i$ is the
flat connection induced by $\al_i$ (using Theorem \ref{fla}) and
$h_i$ is the restriction of $H_i$.\\

Going in the opposite direction, one can complexify a real gerbe to
get a complex gerbe on the same open cover. Through
complexification, a flat real gerbe connection becomes a flat
complex gerbe connection  with 1-connection trivial.

\section{Holonomy bundles on loop spaces and map spaces}\label{hlo}

Consider the free loop space $LX=\{\rho\mid \rho: S^1\to X \;\;
\mr{smooth}\}$. The evaluation map $LX\times S^1\to X$ induces a
homomorphism 
$$H^2(X, G)\to H^2(LX\times S^1, G).$$
 Composing this
with the slant product $H^2(LX\times S^1, G)\to H^1(LX, G)$ over
the generator of $H_1(S^1,G)$, we have the homomorphism
\begin{equation}\label{gr}
\mu: H^2(X, G)\to H^1(LX,G).
\end{equation}
The following results can be viewed as geometric interpretations
of $\mu$ in the cases that $G$ is ${\bf Z}_2, \R^*$, or $\R$.

\begin{thm}\label{geo}
{\em (a)} There is a well-defined line bundle $\wt{l}\to LX$
associated to each gerbe $\mc{G}$ with connection $\na$. The
isomorphism class of $\wt{l}$ is independent of the choice of
$\na$ and the association $[\mc{G}]\mapsto [\wt{l}]$ recovers the
homomorphism $\mu: H^2(X, {\bf Z}_2)\to H^1(LX,{\bf Z}_2)$.

{\em (b)} If $\na$ is flat then $\wt{l}$ carries a natural flat
connection $\wt{\na}$ as well. The holonomies of $\na$ and
$\wt{\na}$ re-establish the homomorphism $\mu: H^2(X, \R^*)\to
H^1(LX,\R^*)$.

{\em (c)} Given a gerbe metric $h$ and a flat connection $\na$ on
$\mc{G}$, each adaptation $\beta$  to $(\na,h)$ corresponds to a
unique metric $\wt{h}$ on $\wt{l}$. Furthermore, the map
$[B]\mapsto [\theta]$ realizes the  homomorphism $\mu:
H^2(X,\R)\to H^1(LX,\R)$, where $B,\theta$ are respectively the
holonomy forms of $(\na,\beta,h)$ and $(\wt{\na},\wt{h})$.
\end{thm}

\nt{\em Proof}. We focus on the constructions, leaving most of the
verifications to the interested reader, since they can be checked
as in the complex gerbe case.

(a) Start with an open cover $\{V_a\}$ of $X$ such that $\{LV_a\}$
covers $LX$, $H^2(V_a,{\bf R}^*)$ is trivial and $V_a\cap V_b$
consists of contractible components for any $a\not=b$. (For
example take a small tubular neighborhood of each loop in $X$. The
cover $\{V_a\}$ does not have to with the cover $\{U_i\}$ where
$\mc{G}$ is locally trivialized.) On each $V_a$, $\mc{G}$
restricts to a trivial gerbe $\mc{G}|_{V_a}$
 because of the
triviality of $H^2(V_a,{\bf Z}_2)$. Applying parts (a), (b) of
Theorem \ref{obj} to $\mc{G}|_{V_a}$ with the gerbe connection
$\na|_{V_a}$, we have an object bundle with  object connection.
Repeat this with $V_b$. On $V_a\cap V_b$, we have now two
restricted objects with connections as well as their  difference
line bundle $\xi_{ab}\to V_a\cap V_b$ with a connection $D_{ab}$.
Introduce a map $g_{ab}:LV_a\cap LV_b\to \R^*$, where at at each
loop $\rho\in LV_a\cap LV_b= L(V_a\cap V_b)$, $g_{ab}(\rho)$ is
the holonomy of $D_{ab}$ along $\rho$. Then our bundle $\wt{l}$ is
defined by the
 transition functions $\{g_{ab}\}$ with respect to the open
cover $\{LV_a\}$.

(b) In this case, the restricted gerbe connection $\na|_{V_a}$ is
flat and has a trivial holonomy from the triviality of
$H^2(V_a,{\bf R}^*)$. By part (c) of Theorem \ref{obj}, the
difference connection $D_{ab}$ is flat, hence  $g_{ab}$ are all
locally constant, because $V_a\cap V_b$ is component-wise
contractible. Thus $\{g_{ab}\}$ gives the expected flat connection
$\wt{\na}$ on $\wt{l}$.

(c) On each $V_a$, apply the proof of part (e) of Theorem
\ref{obj} to the restrictions of $\mc{G}, \na,\beta$, so that we
have a global 1-form $\gamma_a$ (depending on objects). Then the
metric $\wt{h}$ is given by the family of functions $\{\ti{h}_a\}$
on the open cover $\{LV_a\}$, where $\ti{h}_a: LV_a\to \R^+$,
$\rho\mapsto \mr{exp}(2{\int_\rho\gamma_a})$. To see
$\{\ti{h}_a\}$ can be glued via  the transitions $\{g_{ab}\}$, we
just need  to check that at any loop
 $\rho\in LV_a\cap LV_b$,
$\mr{exp}[2{\int_\rho(\gamma_a-\gamma_b)}]=g_{ab}^2(\rho)$, or
equivalently
\begin{equation}\label{exp}
\pm\mr{exp}[{\int_\rho(\gamma_a-\gamma_b)}]= g_{ab}(\rho)
\end{equation}
where the right side  is by definition the holonomy
 of $D_{ab}$ along $\rho$. To compute the left side,
let $\mc{G}$ be locally trivialized on some open cover $\{U_i\}$
of $X$ and let $\{U^a_i\}$ be the induced cover of $V_a$ so that
we have an object, an object connection,
 an object metric $\{l^a_{i}\}, \{\na^a_{i}\}, \{h^a_i\}$ of
$\mc{G}|_{V_a}, \na|_{V_a}, h|_{V_a}$ respectively,
 as in Theorem \ref{obj}.
Then $\gamma_a=\beta_i^a-\theta(\na^a_{i}, h^a_i)$ on
$U_i^a\subset V_a$ from (e) of Theorem \ref{obj}. Likewise,
 $\gamma_b=\beta_{i'}^b-\theta(\na^b_{i'}, h^b_{i'})$ on
$U_{i'}^b\subset V_b$ by working with restrictions to $V_b$. Now
cover the loop $\rho\subset V_a\cap V_b$ with some common
 open sets $W_k=U^a_{i}= U^b_{i'}$. On each $W_k$,
$\gamma_a-\gamma_b=\theta(\na^a_{i}, h^a_i)- \theta(\na^b_{i'},
h^b_{i'})$, since $\beta_i^a=\beta_{i'}^b$, both being
restrictions of the same adaptation $\beta$ from $X$. On the other
hand, $\theta(\na^a_{i}, h^a_i)- \theta(\na^b_{i'},
h^b_{i'})=a_k-\frac{1}{2}d\ln h_k$, where $a_k$ is the connection
1-form of $D_{ab}$
 and $h_k=h^a_i/h^b_{i'}$. Here we use the fact that
$D_{ab}$ is the difference connection of
$\{\na^a_{i}\},\{\na^b_{i'}\}$. Equation (\ref{exp}) now holds
true because with a suitable sign,
$$\pm\mr{exp}[{\int_\rho(\gamma_a-\gamma_b)}]=\pm
\mr{exp}\left[\sum_k\int_{\rho\cap W_k}(a_k-\frac{1}{2}d\ln h_k)
\right]$$
 calculates the holonomy
of $D_{ab}$ along $\rho$, given the (difference) metric
$h_{ab}=\{h_k\}$ on the difference bundle of $\{l^a_{i}\},
\{l^b_{i'}\}$. Indeed one can check that the integrand is the
global holonomy form of $(D_{ab},h_{ab})$ along $\rho$. Compare
with formula (\ref{rh}) in the the proof of Proposition \ref{rel}.

A priori $\gamma_a$ depends  on the choice of various objects
$\{l^a_{i}\}, \{\na^a_{i}\}, \{h^a_i\}$ on $V_a$. But by taking
$a=b$, the above argument shows easily that the function
$\wt{h}_a$ actually does not depend on such choices, since the
difference connection $D_{aa}$ has trivial holonomy. \qed

It is possible to describe the principal bundle of $\wt{l}\to LX$
explicitly, which will be useful later.
\begin{cor}\label{pri}

{\em (a)} The principal $\R^*$-bundle $P_l$ of $\wt{l}$ can be
constructed as follows. At a point $\rho\in LX$, the fiber
$P_{l;\rho}$ consists of equivalence classes of flat object
connections
  of $(\rho^*\mc{G},\rho^*\na)$ on $S^1$, where two flat object
connections  are equivalent if their difference connection has
trivial holonomy.

{\em (b)} If the gerbe connection $\na$ is flat, then the
associated flat connection $\wt{\na}$ on $P_l$ can also be
described explicitly.
\end{cor}

\nt{\em Proof}. (a) For a dimension reason, the pull-back gerbe
connection $\rho^*\na$ on $S^1$ is flat and with trivial holonomy.
By part (c) of Theorem \ref{obj}, the fiber $P_{l;\rho}$ is
well-defined and is acted transitively by the group of isomorphic
flat difference connections on $S^1$. The last group is
$H^1(S^1,\R^*)=\R^*$, hence $P_{l}$ is a principal $\R^*$-bundle
on $LX$. Clearly $P_{l}$ associates with $l$, as the transition
functions of $P_l$ on the open cover $\{LV_a\}$ constructed in the
proof of Theorem \ref{geo} are also given by the holonomy of
difference connections on $S^1$.

(b) Take a path $f:S^1\times[0,1]\to X$ in $LX$ between two loops
$\rho_0,\rho_1\in LX$. Then the pull-back flat gerbe connection
$f^*\na$ has trivial holonomy, since
$H^2(S^1\times[0,1],\R^*)=\{1\}$.
 Part (c) of Theorem \ref{obj} tells us that $(f^*{\mc{G}},f^*\na)$
admits a flat object on $S^1\times[0,1]$. The restricted flat
objects to $\rho_0,\rho_1$ characterize the parallel transport of
$\wt{\na}$ along the path $f$. (Note that here one needs $\na$ to
be flat in order to get a flat $f^*\na$. This is slightly
different from the complex gerbe case, where the pull-back gerbe
connection is automatically flat as the curvature $3$-form has to
vanish on the 2-dimensional manifold $S^1\times[0,1]$.) \qed\\

Next let $\Si$ be an oriented surface with boundary consisting of
$m$ components. On each component $\pa_k\Si$, fix an
orientation-preserving parameterization $S^1\to \pa_k\Si$. Let
$\mr{Map}(\Si,X)$ denote the space of smooth maps from $\Si$ to
$X$. Then we have $m$ natural maps $b_k:\mr{Map}(\Si,X)\to LX$,
where for $f: \Si\to X$, the loop $b_k(f)$ is the composition of
$S^1\to \pa_k\Si$ with $f$. The following extends the main result
in Brylinski \cite{b} to the real case.

\begin{thm}\label{hl}
Suppose $\na=\{\na_{ij}\}$ is  a flat connection  defined on  a
gerbe $\mc{G}=\{l_{ij}\}$.

{\em (a)}  The line bundle
$\ov{l}=\otimes^m_{k=1}b^*_k\wt{l}\to\mr{Map}(\Si,X)$ carries a
canonical trivialization $s$, where $\wt{l}$ is the bundle
constructed in Theorem \ref{geo}.

{\em (b)} The trivialization $s$ is flat with respect to the
pull-back connection $\ov{\na}=\otimes_k b^*_k\wt{\na}$, where
$\wt{\na}$ is the flat connection on $\wt{l}$ induced by $\na$.

{\em (c)} Let  $h=\{h_{ij}\}$ be a metric on $\mc{G}$ and
$\beta=\{\beta_i\}$  an adaptation to $(\na, h)$. Suppose $\wt{h}$
is the metric on $\wt{l}$ induced by $(\na, h, \beta)$. Then $s$
is of norm $1$ with respect to the pull-back metric
$\ov{h}=\otimes_k b^*_k\wt{h}$. Consequently,  $\ov{\na}$ and
$\ov{h}$ are compatible each other.
\end{thm}

\nt{\em Proof}.  For a better presentation, we will work with
principal bundles. Let $P=P_l$ be constructed as in Corollary
\ref{pri} and let
$$P_w=\prod^m_{k=1}b^*_kP=b^*_1P\times\cdots\times b^*_mP$$
be the fiber product bundle on $\mr{Map}(\Si,X)$. $P_w$ has the
structure group $\R^*\times\cdots\times \R^*$ and the associated
bundle $\ov{P}=P_w\times_\eta \R^*$ is exactly the principal frame
 bundle of $\ov{l}$, where $\eta: \R^*\times\cdots\times
\R^*\to \R^*$ is the multiplication homomorphism. (Since both
bundles use the same transition functions.) Thus an element in any
fiber of $\ov{P}$ is represented by some $(d_1,\cdots,d_m) \in
P_w$ inside the equivalence class
\begin{equation}\label{w1}[d_1,\cdots,d_m]=\{(a_1d_1,\cdots,
a_md_m)\in P_w\mid a_1,\cdots,a_m\in\R^*\; \mr{ and } \;\prod_k
a_k=1\}.
\end{equation}

(a) We show $\ov{P}$ is canonically trivial by constructing  a
canonical section. Take any point $f\in\mr{Map}(\Si,X)$, and
consider the pull-back flat gerbe connection $f^*\na$ on $\Si$.
Since $H^2(\Si, \R^*)=\{1\}$, $f^*\na$ must have trivial holonomy
hence  admits a flat object connection $\na_{ob}$ on some object
bundle $\mc{G}_{ob}$ of $f^*\mc{G}$ by Theorem \ref{obj}. Let
$\na^1_{ob}, \cdots, \na^m_{ob}$ be the restrictions to the
boundary components of $\pa{\Si}$. So by Corollary \ref{pri}
 we have an element
$\na^k_{ob}\in P_{f_k}$ in the fiber over  $f_k= b_k(f)\in LX$,
and $\na^k_{ob}\in (b^*_k P)_f$ as well. Now  construct a section
$s$ of $\ov{P}$ by using the equivalence class
$$s(f)=[\na^1_{ob}, \cdots, \na^m_{ob}]\in \ov{P}_f.$$

It remains to show that $s(f)$ is independent of the choice of
objects. Let $\na'_{ob}$ be a second flat object connection of
$f^*\na$ on another object bundle $\mc{G}'_{ob}$ of $f^*\mc{G}$.
These result in the difference flat connection $D$ on the
difference bundle $\xi$ according to Theorem \ref{obj}. The
restricted objects also yield difference flat connections $D_k$ on
the difference bundles $\xi_k$ over all components of $\pa{\Si}$.
Of course $D_k, \xi_k$ are just restrictions of $D, \xi$ to
$\pa_k\Si$. Let $\na_{ob}'^{k}$ denote the restricted flat object
connection to the boundary. By construction of $P_{f_k}$,
$\na_{ob}'^{k}=\al(D_k)\na^k_{ob}$, where $\al(D_k)\in \R^*$ is
the holonomy  of $D_k$ along the loop $f_k$. To show $s(f)$ is
well defined, from (\ref{w1}), we need to check
$\prod^m_{k=1}\al(D_k)=1$. We will apply our holonomy formula
(\ref{rh}) to calculate $\al(D_k)$.

Let $h^t=\{h^t_{ij}\}$ be any gerbe metric on $f^*\mc{G}$. Choose
any object metrics $h_{ob}, h'_{ob}$ on $\mc{G}_{ob},
\mc{G}'_{ob}$, both subordinate to $h^t$. So we have the
difference metric $H$ on $\xi$, with a restriction $H_k$ on
$\xi_k$. By (\ref{rh}), $\al(D_k)=
\pm\mr{exp}\int_{\pa_k\Si}\theta(D_k,H_k)$. Consequently
$$\prod^m_{k=1}\al(D_k)=\left[\prod^m_{k=1}(\pm1)\right][\mr{exp}
\int_{\pa\Si}\theta(D,H)]$$ where we note
$\theta(D,H)|_{\pa_k\Si}=\theta(D_k,H_k)$. Now the total sign
$\prod^m_{k=1}(\pm1)=1$ since the gerbe class $f^*\mc{G}$ is
trivial, and by Stokes' Theorem,
$\int_{\pa\Si}\theta(D,H)=\int_\Si d\theta(D,H)=0$ as the holonomy
form $\theta(D,H)$ is closed. Putting together we arrive at the
desired formula: $\prod^m_{k=1}\al(D_k)=1$.

(b) Recall from Theorem \ref{geo}, the flat connection $\wt{\na}$
on $P$ is constructed via locally constant transition functions
that are holonomy of difference connections of some object
connections. In particular, the constant $t=1$ is a local flat
section of $P$. By using the pull-backs of the afore-mentioned
object connections, it is not hard to see that $t$ pulls back to
$s$ locally. Thus $s$ is locally hence globally flat with respect
to $\ov{\na}$. Alternatively one can use the construction of
$\wt{\na}$ given in the proof of Corollary \ref{pri} to show the
flatness of $s$.

(c) Similar to (a) and (b), the main idea is to use objects on the
whole $\Si$ to calculate the pull-back metrics on boundary
components via restrictions. Choose $\na_{ob}, \mc{G}_{ob}$ as in
(a). Together with the pull-back adaptation $f^*\beta$ and gerbe
metric $f^*h$, we have a global 1-form $\gamma$, as constructed in
the proof of Theorem \ref{obj}. On each $\pa_k\Si$, using the
restrictions of $\na_{ob}, \mc{G}_{ob}, f^*\beta, f^*h$, we have
also the global 1-form $\gamma_k$, which of course is the
restriction of $\gamma$. Then according to Theorem \ref{geo},
under the pull-back metric $\ov{h}$, the norm of $s$ at
$f\in\mr{Map}(\Si,X)$ is
$$\prod_k\mr{exp}(2\int_{\pa_k\Si}\gamma_k)=
\mr{exp}(2\int_{\pa\Si}\gamma).$$
 By Stokes' Theorem, $\int_{\pa\Si}\gamma=\int_\Si d\gamma=\int_\Si B$,
where $B$ is the holonomy form of $(f^*\mc{G},f^*\na,f^*\beta)$.
Since the holonomy class of $f^*\na$ is trivial, $\int_\Si B=0$ by
Theorem \ref{cur}, and the norm $|s(f)|=1$ consequently. \qed

\begin{rem}
To be logically correct, we check that the argument in (a) is
independent of the gerbe metric $h^t$ on $f^*\mc{G}$. Let
$h^{tt}=\{h^{tt}_{ij}\}$ be a second choice. Then there exist
positive functions $f_{ij}$ such that
$h^{tt}_{ij}=f_{ij}h^{t}_{ij}$. Since the sheaf $\un{\R}^+$ of
positive functions is fine, there is a Cech 0-cocycle $\{f_i\}$
with coboundary $\{f_{ij}\}$. Multiplying $h_{ob}, h'_{ob}$ by
$\{f_i\}$ we get two object metrics subordinate to $h^{tt}$.
However their difference metric on $\xi$  remains the same, and so
does the rest argument. We can also check that the argument does
not depend upon the choice of object metrics $h_{ob}, h'_{ob}$. On
$\mc{G}_{ob}, \mc{G}'_{ob}$ take another pair of object metrics,
still  subordinate to $h^t$. For the resulting  metric $H^\circ$
on $\xi$, the 1-form
 $\theta(D,H^\circ)$  differs from $\theta(D,H)$ by an exact 1-form
on $\Si$. Hence the integral on $\pa\Si$ remains the same.
\end{rem}

Let us now incorporate the real picture of the section. Suppose
$\si:M\to M$ is a smooth involution with fixed point set
$M_{\R}=X$. For a loop $\rho\in LM$, define $\ov{\rho}\in LM$ by
$\ov{\rho}(t)=\si(\rho(t))$. Then the fixed loop space
$(LM)_{\R}=LX$. Given a complex gerbe $\mc{G}^c$ on $M$ with
connection $\na^c=\{\na^c_{ij}, F_i\}$, a well-defined complex
line bundle $\wt{L}\to LM$ is constructed in \cite{b, c} together
with a connection $\wt{\na}^c$. The following sums up the relation
between the two gerbe pictures on $M$ and $X$. Its proof is
 essentially book-keeping.

\begin{pro}\label{lx}
 {\em(a)} Suppose the gerbe $\mc{G}^c$ and  connection
$\na^c$ are both Real. Then the associated complex line bundle
$\wt{L}$ and connection $\wt{\na}^c$ are also Real. Taking real
parts, we have a real line bundle $\wt{L}_{\R}\to LX$ and a
connection $\wt{\na}^c_{\R}$. They are naturally identified with
the bundle $\wt{l}$ and connection $\wt{\na}$, which are
constructed by using the real gerbe $\mc{G}$ and connection $\na$.
If $\na^c$ is flat or unitary,  then  $\wt{\na}^c_{\R}=\wt{\na}$
is flat as well.

{\em(b)} Each Hermitian gerbe metric $H$ on $\mc{G}^c$ induces a
unique metric $\wt{h}$ on $\wt{l}$.

{\em(c)} Suppose further that $\si$ is anti-holomorphic on a
complex manifold $M$ and $\mc{G}^c$ is a holomorphic gerbe with a
holomorphic structure $\al$. Then we have a second line bundle
$\wt{l}^\al\to LX$ with a flat connection $\wt{\na}^\al$. Moreover
$\wt{l}^\al$ is isomorphic with $\wt{l}$.
\end{pro}

\nt{\em Proof}. (a) By construction in \cite{b,c}, $\wt{L},
\wt{\na}^c$ are obviously Real. To see $\wt{L}_\R=\wt{l}$ amounts
to comparing the holonomies of a Real connection and its real
part.

Let $G=dF_i$ denote the curvature 3-form of $\na^c$. According to
\cite{b, c} the curvature $K$ of $\wt{\na}^c$ is evaluated at
$\rho\in LM$  as
$$K(u,v)=\mr{exp}\int_\rho i_{u,v}G,$$
where $u,v\in T_\rho(LM)$ are two vector fields along $\rho$ and
$i_{u,v}$ is the contraction. If $\na^c$ is flat, $G=0$, then
$K=0$ and $\wt{\na}^c$ is flat,  and so is its real part
$\wt{\na}^c_\R$. (But note that the gerbe connection $\na$ on
$\mc{G}$ may not be flat.)

If $\na^c$ is unitary, then $\na$ is flat and so is
$\wt{\na}=\wt{\na}^c_\R$.

(b) $H$ restricts to a Riemann metric $h=\{h_{ij}\}$ on $\mc{G}$.
Take the unique compatible connection $\na^h_{ij}$ of $h_{ij}$ on
$l_{ij}$ with $\theta(\na^h_{ij},h_{ij})=0$. Together with the
trivial adaptation $\beta_i=0$, we get a metric $\wt{h}$ by
Theorem \ref{geo}. Note that there is no Hermitian metric on
$\wt{L}$ directly induced
 by $H$.

(c) Let $\na^\al$ be the flat gerbe connection on $\mc{G}$ that is
associated to $\al$. Then using $(\mc{G},\na^\al)$, we have a real
line bundle $\wt{l}^\al\to LX$ with flat connection $\wt{\na}^\al$
by Theorem \ref{geo}. One sees that $\wt{l}^\al,\wt{l}$ have the
same isomorphism type, since both are determined by that of
$\mc{G}$. \qed

Note that for the sake of our next discussion, we have extended
$\si:M\to M$ to $\bar{\si}: LM\to LM$ by using the identity map on
$S^1$. One could also use the antipodal or complex conjugation
maps on $S^1$, in which cases we no longer have the fixed loop
spaces $(LM)_\R=LX$. Such cases also appear interesting to study.

Our next purpose is to study the Real analogue of Theorem \ref{hl}
in connection with map spaces. Let $\Si$ be a closed complex curve
with a real structure (namely an anti-holomorphic involution). The
real part $\Si_{\R}$ consist of $m$ circles $\Si_1,\cdots,\Si_m$.
Continue assuming $\si: M\to M$ to be a smooth involution such
that $X=M_\R$ is a smooth manifold, as before. Fix any orientation
on each $\Si_k$ and orientation-preserving diffeomorphism $S^1\to
\Si_k$. Obviously $\mr{Map}(\Si,M)$ inherits a natural smooth
involution and the real part $\mr{Map}_{\R}(\Si,M)$ contains all
the Real maps from $\Si$ to $M$. Using the fixed diffeomorphism
$S^1\to \Si_k$ we have a map $b_k:\mr{Map}_{\R}(\Si,M)\to LX$.
Take a Real complex
 gerbe $\mc{G}^c$ with connection $\na^c$ on $M$ and
let $\wt{l}\to LX$ be the real line bundle from Proposition
\ref{lx}. Set the real line bundle $\ov{l}\to\mr{Map}_{\R}(\Si,M)$
to be $\otimes_k b_k^*\wt{l}$. The main question here is about the
triviality of $\ov{l}$. (Incidentally, reversing orientation on
any $\Si_k$ will not change the isomorphism type of $\ov{l}$.)
Certainly the answer does not come directly from Theorem \ref{hl}:
for one thing $\mr{Map}_{\R}(\Si,M)$ explicitly involves $M$ while
$LX$ does not. In fact we are able to answer the question only in
a special case.

Recall that the complement $\Si\backslash \Si_\R$ contains at most
two connected components; cf. Wilson \cite{wi} for example. The
real structure on $\Si$ is called dividing or non-dividing,
depending on whether there are two or one components. In the
dividing case, there is a well-defined pair of opposite
orientations on all $\Si_\R$, which are induced by orientations on
the two components.

\begin{thm}\label{rhl} Let $\Si$ be a closed complex curve
carrying a dividing real structure and $M$ carrying a smooth
involution with real part $X$. Use the real part $\Si_\R$ to
define  maps $b_k: \mr{Map}_{\R}(\Si,M)\to LX$ as above.

{\em (a)} Suppose $\wt{l}\to LX$ with connection $\wt{\na}$ is
constructed as in part (a) of Proposition \ref{lx}. Then
$\ov{l}:=\otimes_k b_k^*\wt{l}$ is canonically trivialized. If the
gerbe connection $\na^c$ is flat or unitary, then the
trivialization is flat with respect to the pull-back connection
$\ov{\na}=\otimes_k b^*_k\wt{\na}$.

{\em (b)} Suppose we are
in the complex set up of part (c) of Proposition \ref{lx}, so we
have the real line bundle $\wt{l}^\al\to LX$ with flat connection
$\wt{\na}^\al$. Then  the pull-back bundle $\ov{l}^\al=\otimes_k
b_k^*\wt{l}^\al$ has a flat trivialization with respect to the
connection $\ov{\na}^\al=\otimes_k b^*_k\wt{\na}^\al$.
\end{thm}

\nt{\it Proof}. (a) Pick one of the two components  in
$\Si\backslash \Si_\R$ and orient $\Si_\R$ accordingly. Fix any
orientation-preserving diffeomorphisms $S^1\to \Si_k$ where
$\Si_\R=\coprod^m_{k=1}\Si_k$. Label the chosen component as
$\Si^+$, with the boundary included. Extending the usual complex
gerbe case for surface with boundary, we consider the map
$b_k:\mr{Map}(\Si,M)\to LM$ and the complex line bundle
$\ov{L}=\otimes b^*_k\wt{L}$, where $\wt{L}$ is as in Proposition
\ref{lx}. (This is a generalization of the usual complex case, as
$\Si$ does not have a boundary -- the real part $\Si_\R$ plays the
same role here.) Since $\Si_k$ contains real points only, the map
$b_k$ is clearly Real with respect to the natural real structures.
The restriction to real parts $b_k:\mr{Map}_\R(\Si,M)\to LX$ is
exactly the one defined early. Note $\ov{L}$ has a real lifting as
$\wt{L}$ is Real, and $\ov{l}\to \mr{Map}_\R(\Si,M)$ is exactly
the real part of $\ov{L}$ under such real lifting.

To show $\ov{l}$ is canonically trivial, it is enough to show
$\ov{L}$ is so with a Real trivialization. This in turn follows
essentially the same kind of proof as  Theorem \ref{hl}. Indeed,
at $f\in\mr{Map}(\Si,M)$, let $f_+$ be its restriction to $\Si^+$.
For dimension reason, the pull-back gerbe connection $f^*_+\na^c$
is flat. And its holonomy is trivial as well as $H^3(\Si^+,{\bf
C}^*)$ is certainly trivial. Thus $f^*_+\mc{G}^c$ admits object
bundle $\mc{G}^c_{ob}$ with connection $\na^c_{ob}$ by the complex
version of Theorem \ref{obj} (see \cite{b, c}). The restrictions
$\na^{c,1}_{ob},\cdots,\na^{c,m}_{ob}$ of $\na^c$ to the circles
$\Si_1,\cdots,\Si_m$ yield an element $s(f)$ in the fiber of the
principal bundle $\ov{P}^c$ associated to $\ov{L}$. That $s(f)$ is
independent of the choice of objects follows again from the
holonomy formula along circles and the Stokes Theorem applied to
$\Si^+$ with boundary $\Si_\R$. Consequently we have $\ov{P}^c$
hence $\ov{L}$ canonically trivialized by the section $s$. Since
$\mc{G}^c,\na^c$ are both Real, it is not hard to check that $s$
is invariant under the real lifting on $\ov{L}$: essentially this
is due to the fact that the other component $\Si^-$ of
$\Si\backslash \Si_\R$ induces the opposite orientation on
$\Si_\R$ as $\Si^+$ does. Restricting to the real
part$\mr{Map}_\R(\Si, M)$, $s$ becomes the canonical
trivialization of $\ov{l}$.

When $\na^c$ is flat or unitary, the connection $\wt{\na}$ is flat
hence $\ov{\na}$ is flat as well. The trivialization $s$ is
$\ov{\na}$-flat, following a similar proof to Theorem \ref{hl}.

(b) By Proposition \ref{lx}, $\ov{l}^\al$ is isomorphic with
$\ov{l}$ hence is trivial. Moreover the induced trivialization
from $\ov{l}$ is $\ov{\na}^\al$-flat. \qed

When the real structure on $\Si$ is non-dividing, we are unable to
obtain any definitive result but conjecture that $\ov{l}$ is not
trivial in general. A basic result in real algebraic geometry says
that $\Si$ is dividing iff the class $[\Si_\R]\in H_1(\Si, {\bf
Z}_2)$ is trivial, see \cite{wi} for example. Thus in the
non-dividing case, the Poincare dual $PD[\Si_\R]\in H^1(\Si, {\bf
Z}_2)$ associates with a non-trivial real line bundle on $\Si$. We
conjecture that $\ov{l}$ should be related to this line bundle.\\

We end the paper with a speculation concerning real Gromov-Witten
invariants. Continue using the real complex set-up from part (b)
of Theorem \ref{rhl}. Let $\mc{M}_\R\subset\mr{Map}_\R(\Si,M)$ be
the moduli space of real holomorphic curves. Take a multi-degree
cohomology class $\beta=(\beta_1,\cdots,\beta_m)$, where
$\beta_k\in H^{n_k}(LX,\wt{l})$ and $\wt{l}=\wt{l}^\al$ carries
the flat connection $\wt{\na}^\al$ so that the cohomology with
twisted coefficients is defined. Assume the total degree of
$\beta$ is $\sum_k n_k=\dim \mc{M}_\R$. Set
$b^*\beta=b^*_1\beta_1\cup\cdots\cup b_m^*\beta_m$, which is an
ordinary cohomology class of degree $\dim \mc{M}_\R$ on
$\mr{Map}_\R(\Si,M)$, since the pull-back bundle $\ov{l}$ is
trivial. Then one might attempt to define a real type of
Gromov-Witten invariant with respect to a real holomorphic gerbe
$(\mc{G}^c,\al)$ as the map
$$GW_{\R,\al}:  H^{n_1}(LX,\wt{l})\times\cdots\times H^{n_m}(LX,\wt{l})
\longrightarrow \R$$ by sending $\beta$ to
$(b^*\beta,[\ov{\mc{M}}_\R])$, where $\ov{\mc{M}}_\R$ is a
compactified real moduli space. Here implicitly we have asserted
that the real moduli space $\mc{M}_\R$ is orientable, which seems
plausible under our assumption that $\Si$ has a  dividing real
structure (as in Theorem \ref{rhl}). Compare with Katz-Liu \cite{kl}.

In the non-dividing case, $b^*\beta$ belongs to the cohomology of
twisted coefficients $\ov{l}$ and $\mc{M}_\R$ may not be
orientable either. However, if the orientation bundle of $\mc{M}_\R$
matches $\ov{l}$, then the above definition for $GW_{\R,\al}$ can
make sense again.

A Gromov-Witten invariant with respect to a complex gerbe has been
introduced by Pan, Ruan, and  Yin \cite{pry}. We hope to return to
the study of the real case in a future work.

\vspace{10mm}

\nt{\bf Acknowledgment.} The author wishes to thank Yongbin Ruan
for introducing him to the topic of gerbes.

\end{document}